\documentclass[11pt, a4paper, oneside]{amsart}
\usepackage[utf8]{inputenc}
\usepackage{amsmath}
\usepackage{float}
\usepackage{amssymb}
\usepackage{amsthm}
\usepackage{tkz-euclide}
\usepackage{pgfplots}

\pgfplotsset{compat=1.10}
\usepgfplotslibrary{fillbetween}
\usepackage{mathtools}
\usepackage{setspace}
\usepackage[a4paper, total={6in, 8in}]{geometry}
\newcommand*\Laplace{\mathop{}\!\mathbin\bigtriangleup}

\usepackage{enumitem}   
\newtheorem{theorem}{Theorem}
\newtheorem{proposition}[theorem]{Proposition}
\newtheorem{corollary}[theorem]{Corollary}
\newtheorem{lemma}[theorem]{Lemma}
\theoremstyle{definition}
\newtheorem{definition}{Definition}
\theoremstyle{remark}

\usepackage{tikz}
\usepackage{bm}
\usepackage{caption}
\usepackage{subcaption}
\usepackage{graphicx}
\graphicspath{ {./pictures/} }
\usepackage{import}
\usepackage{pdfpages}
\usepackage{xifthen}
\usepackage{indentfirst}
\usepackage{transparent}
\DeclareMathOperator{\ran}{ran}

\usepackage{enumitem}
\usepackage{lipsum}
\usepackage[colorlinks, citecolor = blue]{hyperref}
\usepackage{cleveref}
\title[Gravity water waves over constant vorticity flows]{Gravity water waves over constant vorticity flows: from laminar flows to touching waves}
\author{Francisco Gonçalves}
\newlist{menumerate}{enumerate}{2}
\setlist[menumerate,1]{label=(\alph*)}
\setlist[menumerate,2]{label=(\roman*)}
\begin{document}
\begin{abstract}  
In a recent paper, Hur \& Wheeler \cite{paper2} proved the existence of periodic steady water waves over an infinitely deep, two-dimensional and constant vorticity flow under the influence of gravity. These solutions include overhanging wave profiles, some of which exhibit surfaces that touch at a point and thereby enclose a bubble of air. We extend these results by formulating a problem that encompasses both infinitely deep and finitely deep flows, and by proving the existence of a continuous curve of water waves that connects a laminar flow to a touching wave for fixed, nonzero gravity. This implies the existence of a wave profile featuring a vertical tangent at a point, which is not overhanging, and is referred to as a breaking wave. We also study the behaviour of critical layers, which are points where the horizontal velocity vanishes, near the surface. In particular, this result holds for arbitrary vorticity.
\end{abstract}
\maketitle
\footnotetext[1]{This work is based on my master's thesis. I would like to express my sincere gratitude to my supervisors Evgeniy Lokharu and Erik Wahlén for their invaluable guidance, feedback, and support throughout the completion of the master's thesis. I would also like to thank Anna Geyer for a careful reading of the introduction. This work was partially supported by grant VI.Vidi.223.019 from the Dutch Research Council (NWO).}
\section{Introduction}

\par We study periodic water waves travelling over a two-dimensional, incompressible, inviscid flow with constant vorticity. The effect of gravity is included, while surface tension is neglected. We prove the existence of waves with an overhanging profile for flows of finite depth, including profiles whose surfaces touch at a point, referred to as touching waves. Furthermore, for both finite and infinite depth, and for fixed nonzero gravity, we establish the existence of a continuous curve of solutions connecting a laminar flow to a touching wave.

\par In 1957, Crapper \cite{crapper}  discovered an explicit and smooth curve of exact solutions for capillary waves propagating in an irrotational flow in the absence of gravity. As the amplitude increases, the crests of Crapper's solutions become flatter while the troughs become sharper. This behavior contrasts with that of gravity waves. The endpoint of Crapper's curve is a touching wave. Akers, Ambrose, and Wright \cite{akers2014gravity}, as well as de Boeck \cite{deboeck2014existence}, constructed small-amplitude gravity waves close to Crapper’s solutions. Córdoba, Enciso, and Grubic \cite{CordobaEncisoGrubic} further extended this work by constructing a touching wave.

\par In the setting of constant vorticity with zero surface tension and zero gravity, numerical studies by Dyachenko and Hur \cite{dyachenko_hur_2019,https://doi.org/10.1111/sapm.12250} and Hur and Vanden-Broeck \cite{HUR2020190} found evidence of explicit solutions with the same wave profile as Crapper’s solutions. Hur and Wheeler \cite{paper1} later established the existence of such solutions analytically. More recently, Crowdy \cite{crowdy_2023} showed that these solutions can be understood as part of a broader framework of explicit solutions in the constant vorticity case under the same assumptions.

\par Using an Implicit Function Theorem argument, Hur \& Wheeler \cite{paper2} constructed overhanging and touching waves for small gravitational effects in a constant-vorticity flow. In this paper, we take these matters further, by considering both finitely and infinitely deep flows, and by applying local bifurcation theory and compactness arguments to improve upon their results. This leads to a smooth curve of solutions from the laminar flow up to a touching wave for fixed small nonzero gravity. Along this curve, we obtain waves whose profiles are vertical at a point but never overhanging, referred to as breaking waves.

\par We also analyse the local behaviour of critical layers (points where the horizontal velocity of the fluid vanishes) at the surface of breaking and overhanging waves and how they interact with gravity. In particular, this result holds for arbitrary vorticity.

\par We begin by introducing and reformulating the water wave problem in Section \ref{Section: prelim and reform}. This is followed by an extension of the result of Hur and Wheeler to the finite depth case, presented in Section \ref{section: overhanging waves}. Next, in Section \ref{section: cont curve}, we construct a continuous curve of solutions starting from the laminar flow and continuing to the touching wave. Finally, Section \ref{section: crit layers} is devoted to analysing the local behaviour of critical layers at points where the wave has a vertical tangent.

\section{Preliminaries and reformulation}\label{Section: prelim and reform}
\par We consider the water wave problem for a periodic, steady wave travelling over a two-dimensional, incompressible, inviscid fluid under the influence of gravity while neglecting surface tension. For simplicity, we assume the fluid has unit density. In Cartesian coordinates, we assume the $x$-axis points in the direction of wave propagation, while gravity acts in the negative $y$-direction. We consider a moving frame of reference so that the flow is stationary. Suppose the flow occupies a region $D$ in the $x,y$ plane and is bounded from above by a free surface $S$. Using a stream function, $\psi$, we can express this problem as
\begin{subequations}
\label{conditions pre}
\begin{align}
&\Laplace\psi=-\omega &\text{in }D,\label{conditions poisson}\\
&\psi=0  &\text{on }S,\label{kinematic conditions1}\\
&\frac{1}{2}|\nabla\psi|^2+ g y= b &\text{on }S,\label{dynamic boundary condition in conditions1}
\end{align}
\end{subequations}
 where $\omega$ denotes the vorticity, which is assumed to be constant in $D$. Moreover, $g$ denotes the gravitational constant and $b$ is used to denote a Bernoulli constant. Note that $S$ is a free boundary, \eqref{kinematic conditions1} and \eqref{dynamic boundary condition in conditions1} are, respectively, the kinematic and dynamic boundary conditions. We assume that both $D$ and $\psi$ are $\tfrac{2\pi}{k}$ periodic in the $x$ direction and symmetric with respect to the vertical lines beneath the wave crest and trough. In addition, we consider two alternative boundary conditions at the bottom:
\begin{subequations}
\label{bottom conditions pre}
\begin{align}
 &\nabla\psi-(0,-\omega y-c)\rightarrow (0,0) &\text{as }y\rightarrow -\infty,\label{kinematic at the infinity bottom pre}\\
 &\psi=q &\text{on }y=-h,\label{kinematic at the finite bottom pre}
 \end{align}
\end{subequations}
 where $c$ denotes the wave speed and $q$ is an a priori unknown constant. Note that in the finite depth case, the wave is also assumed to travel with speed $c$. Equation (\ref{kinematic at the infinity bottom pre}) corresponds to the infinite depth case, where $D$ is unbounded from below and (\ref{kinematic at the finite bottom pre}) corresponds to the finite depth case, where $D$ is bounded from below by $\{y=-h\}$, where $h$ is a positive constant. In fact, in the finite depth case, the problem is invariant under a vertical translation after a suitable modification of $B$.
 \par We now make the change of variables
\begin{equation*}
x\rightarrow kx, \qquad\qquad y\rightarrow ky,\qquad\qquad \psi\rightarrow \frac{k}{c}\psi
\end{equation*}
 and introduce the dimensionless parameters
\begin{equation*}
\Omega=\frac{\omega}{ck},\qquad G=\frac{g}{kc^2},\qquad B=\frac{b}{c^2 },  \qquad H=-hk,\qquad Q=\frac{qc}{k}.
\end{equation*}
 In these non-dimensional variables, (\ref{conditions pre}) becomes
\begin{subequations}
\label{conditions}
\begin{align}
&\Laplace\psi=-\Omega &\text{in }D,\\
&\psi=0  &\text{on }S,\label{kinematic conditions}\\
&\frac{1}{2}|\nabla\psi|^2+ G y= B &\text{on }S\label{dynamic boundary condition in conditions}
\end{align}
\end{subequations}
and the bottom boundary conditions (\ref{bottom conditions pre}) become
\begin{subequations}
\label{conditions bottom}
\begin{align}
&\nabla\psi-(0,-\Omega y-1)\rightarrow (0,0) &\text{as }y\rightarrow -\infty,\label{kinematic at the infinity bottom}\\
&\psi=Q &\text{on }y=-H.\label{kinematic at the finite bottom}
\end{align}
\end{subequations}
 This completes the non-dimensional formulation of the problem. In the next section, we reformulate it using tools from complex analysis.
\subsection{Reformulation}
For any $d>0$ possibly infinity, we define
\begin{equation*}
\mathbb{S}_d:= \{\alpha+i\beta \in\mathbb{C};0>\beta>-d\}.
\end{equation*}
\begin{lemma}\label{conformal map lemma}
Consider a $2\pi$ periodic domain $D$ that is bounded from above by $S$ and is either unbounded from below or bounded by $\{y=-H\}$. Then there exists a unique $d$ and a conformal map $z$ from $\mathbb{S}_d$ onto $D$ satisfying the following conditions.
\renewcommand{\theenumi}{(\roman{enumi})}
\begin{enumerate}
    \item $z$ extends continuously to $\{\beta=0\}$ and maps $\{\beta=0\}$ to $S$;
    \item $x(\alpha+i\beta)-\alpha$ and $y(\alpha+i\beta)$ are $2\pi$ periodic functions of $\alpha$;
    \item if $H<\infty$ then $d<\infty$ and the continuous extension of $z$ maps $\mathbb{S}_d\cup\{\beta=-d\}$ to $D\cup \{y=-H\}$ continuously;
    \item if $H=\infty$ then $d=\infty$ and $z(\alpha+i\beta)-(\alpha+i\beta)\rightarrow 0 \text{  as  } \beta\rightarrow -\infty.$
\end{enumerate}
 where $x(\alpha+i\beta)$ and $y(\alpha+i\beta)$ denote the real and imaginary parts of $z(\alpha+i\beta)$, respectively.
\end{lemma}
\par For a proof of the finite-depth case, see the appendix of \cite{conformalmapConstantinVarvaruca}. The proof for the infinite-depth case is very similar. 

We now state a result concerning the analyticity of the free surface $S$ and the analytic extension of the stream function $\psi$. Constantin and Escher \cite{analyticoriginal} established a primary result on the analyticity of the free surface under the assumptions that the fluid contains no stagnation points and that the surface can be represented as the graph of a function. Notably, their analysis permits the vorticity to be a real analytic function. Later, Aasen and Varholm \cite[Theorem 2.5]{analyticsecond} improved this result for the affine vorticity case, showing that the no-stagnation assumption is required only on the boundary. Our proof follows the strategy of \cite{analyticsecond} and uses a local argument that accommodates free surfaces that are not representable as single-valued graphs.
\begin{proposition}\label{analytic proposition}
Let $z$ be the conformal map obtained in Lemma \ref{conformal map lemma} and let $\psi$ be a solution of (\ref{conditions}). Assume that the free surface $S$ is a $C^1$ curve and that $z_\alpha$ has no zeros on the set $\{\beta=0\}$. Then the following three properties hold.
\renewcommand{\theenumi}{(\roman{enumi})}
\begin{enumerate}
    \item The free surface $S$ is real analytic;
    \item The map $z$ has a complex analytic extension to an open neighbourhood of $\mathbb{S}_d\cup \{\beta=0\}$;
    \item The stream function $\psi$ has a real analytic extension that satisfies $\Laplace\psi=-\Omega$ to an open neighbourhood of $\mathbb{S}_d\cup \{\beta=0\}$.
\end{enumerate}
\end{proposition}

\begin{proof}\hfill
\renewcommand{\theenumi}{(\roman{enumi})}
\begin{enumerate}
    \item Let $p$ be an arbitrary point on $S$. There exists an open set $U$ contained in $D$, such that $p\in \partial U \cap S$ and $S\cap U$ can be represented either as $\left(x,\eta(x)\right)$ or $\left(\eta(y),y\right)$ for some function $\eta$. Applying \cite[Theorem 2]{nirenbergregularity} with $u$ as the stream function and $g$ as the dynamic boundary condition. Note here the assumption of $z_\alpha\neq 0$ on the boundary is equivalent to $\nabla\psi\neq 0$ on the boundary. This yields that $\partial U \cap S$ is analytic. Since $p$ was arbitrary, it follows that $S$ is real analytic. 

    \item This follows from applying \cite[Theorem 2.2, page 299]{lang1985complex}.
    \item Since $\Laplace \psi +\Omega=0$ in $D$ and $\psi=0$ on the real analytic curve $S$, we can apply \cite[Theorem A]{Morrey1957OnTA} to yield the claim.
\end{enumerate}
\end{proof}
Proposition \ref{analytic proposition} ensures that both the free surface and the associated stream function extend analytically to nearby regions, this is necessary for studying the local behaviour of critical layers in Section \ref{section: crit layers}.

We use $\mathcal{H}_\infty$, to denote the standard periodic Hilbert transform. Next, we introduce the periodic Hilbert transform on a strip.
\begin{definition}[Periodic Hilbert transform on a strip]
For $d<\infty$, we define the periodic Hilbert transform on $\mathbb{S}_d$, $\mathcal{H}_d:L^2_{2\pi}\rightarrow L^2_{2\pi}$ by its action on the orthonormal basis $\{\exp(inx)\}_{n\in\mathbb{Z}}$, as follows
\begin{equation*}
\mathcal{H}_d(\exp(inx))=\begin{cases}-i\coth(nd)\exp(inx)&\text{if }n\neq 0,\\
0&\text{if }n= 0,
\end{cases}
\end{equation*}
 where $\coth$ denotes the hyperbolic cotangent.
\end{definition}
Using the properties of the Hilbert transform and Lemma \ref{conformal map lemma}, we can reformulate (\ref{conditions}) and (\ref{conditions bottom}) as
\begin{equation}\label{finite depth babenko}
\frac{1}{2}(1+\Omega (y+y\mathcal{H}_dy_\alpha-\mathcal{H}_d(yy_\alpha)))^2=(B-Gy)((1+\mathcal{H}_dy_\alpha)^2+y_\alpha^2)\qquad \text{on }\beta=0.
\end{equation}
We refer to \cite{https://doi.org/10.1111/sapm.12250} for a thorough derivation of (\ref{finite depth babenko}). 

A solution of \eqref{finite depth babenko} gives rise to a solution of the physical problem provided that
\begin{equation*}
z(\alpha+i0)=\alpha+(i+\mathcal{H}_d)(y)(\alpha+i0)
\end{equation*}
is injective and $z_\alpha(\alpha+i0)\neq 0$ for all $\alpha \in\mathbb R$. In the finite depth case $H$ is determined via the inversion of the map $z$, this $H$ is unique up to a shift of $D$ in the vertical direction.

We introduce the variable $\zeta=\exp(-i(\alpha+i\beta))$ and the domains $\mathbb{D}:=\{\zeta\in\mathbb{C}:|\zeta|\leq 1\}$ and 
\begin{equation*}
\mathbb{A}_d:=
\begin{cases}
\{\zeta\in\mathbb{D}:|\zeta|>\exp(-d)\}&\text{if }d\neq \infty,\\
\mathbb{D}&\text{if } d=\infty
\end{cases}
\end{equation*}
and define $w:\mathbb{A}_d\rightarrow \mathbb{C}$ to be such that
\begin{equation}\label{wdecomptwo}
z(\zeta)=i\log\zeta + w(\zeta).
\end{equation} 
This decomposition transforms (\ref{finite depth babenko}) into
\begin{equation*}
\left(B-G\Im(w)\right)=\frac{1}{2}\frac{\left(1+\Omega(\Im w+\mathcal{Q}^d(w)\right)^2}{|1-i\zeta w_\zeta|^2}\qquad \text{for }|\zeta|=1,
\end{equation*}
where $\mathcal{Q}^d$ consists of a commutator term, namely 
\begin{equation*}
\mathcal{Q}^d(w(\zeta)):= \Im(w)\mathcal{H}_d\left(\Im(-i\zeta w_\zeta)\right)-\mathcal{H}_d\left(\Im w\Im(-i\zeta w_\zeta)\right)\qquad \text{for }|\zeta|=1.
\end{equation*}

\subsection{Operator formulation}
\par We define the Banach spaces
\begin{align*}\label{spaces fin}
X&:= \{w\in C^{3+a}(\partial\mathbb{D},\mathbb{C}):w(\overline{\zeta})=w(\zeta)\},\\
Y&:= \{w\in C^{2+a}(\partial\mathbb{D},\mathbb{R}):w(\overline{\zeta})=w(\zeta)\}.
\end{align*}
For a given $w\in X$ and constant $d$, we can extend $w$ to a function $\Tilde{w}$ that is harmonic on $\mathbb{A}_d$ and satisfies the following conditions
\renewcommand{\theenumi}{(\roman{enumi})}
\begin{enumerate}
    \item $\Tilde{w}=w$ on $\partial \mathbb{D}$;
    \item If $d\neq \infty$, $\Tilde{w}$ is constant on the inner boundary of $\mathbb{A}_d$. Moreover, the constant is chosen so that the mean over the inner circle of $\mathbb{A}_d$ is equal to the mean over the outer circle of $\mathbb{A}_d$.
\end{enumerate}
 We can then use $\Tilde{w}$ and its harmonic conjugate to define a holomorphic function on $\mathbb{A}_d$. Next, we introduce the map $E:X\times (0,\infty]\rightarrow Y$ which assigns to each pair $(w,d)$ the function that is holomorphic on $\mathbb{A}_d$ and whose imaginary part satisfies conditions (i) and (ii) when $d\in(0,\infty)$; in the case $d=\infty$. $E(w,\infty)$ on $D$ with its imaginary part satisfying condition (i).
 We then define the subset
\begin{align*}
U&:=\left\{w\in X: 1-i\zeta \frac{\partial E(w,\infty)}{\partial \zeta}\neq 0 \right\}
\end{align*}
and the operator $\mathcal{G}:U\times\mathbb{R}\times (-\infty,\tfrac{1}{3})\times \mathbb{R}\rightarrow Y$
\begin{equation}\label{more general of operators}
\mathcal{G}(w,G,a,l)= 
\frac{1}{2}\frac{\left(1+\Omega(a)\left(w+\mathcal{Q}^\frac{1}{l^2}\left(E(w,\frac{1}{l^2})\right)\right)\right)^2}{|1-i\zeta  \partial_\zeta E(w,\frac{1}{l^2})|^2}-\left(B(a)-G  w\right) \qquad\text{for }|\zeta|=1,
\end{equation}
 where
\begin{equation*}
\mathcal{Q}^\frac{1}{l^2}=
\begin{cases}
\mathcal{Q}^\frac{1}{l^2}&\text{if }l\neq 0,\\
\mathcal{Q}^\infty &\text{if }l= 0
\end{cases}
\end{equation*} 
and
\begin{equation}\label{parametersnograv}
    \Omega(a):= \frac{1-a}{1-3a},\qquad B(a):=\frac{1}{2}\left(\frac{1+a}{1-3a}\right)^2.
\end{equation}
A similar operator was studied by Hur \& Wheeler \cite{paper2}. The difference between our operator and theirs is the additional parameter $l$.

\section{Overhanging waves over a finitely deep flow}\label{section: overhanging waves}
We define
\begin{equation}\label{sol of F 0 grav}
w(a)(\zeta):=-\frac{4i \sqrt{a}\zeta}{1+\sqrt{a} \zeta}\qquad \text{for }|\zeta|=1,
\end{equation}
while noting this is the form Hur and Wheeler's exact solutions \cite{paper1} take in our notation. These solutions are overhanging for $a\in(a_{\text{crit}},a_{\text{max}})$. Here $a_{\text{crit}}=(\sqrt{2}-1)^2$ can be explicitly computed by solving $x_\alpha(\alpha,0;a):=\Re(i\log \zeta+w(a)(\zeta)=0$ and $a_{\text{max}}\approx \sqrt{0.454670016452010}$ was given in \cite{paper2}.

It follows from \cite{paper2} that $\mathcal{G}(w(a),0,a,0)=0$ for all $a\in(0,\frac{1}{4})$. The partial Fréchet derivative of $\mathcal{G}$ can be computed to be
\begin{equation*}
\begin{aligned}
D_w\mathcal{G}[w,G,a,0]v=&\frac{1+\Omega(a)(w+\mathcal{Q}(E(w,\infty))}{|1-i\zeta \partial_\zeta E(w,\infty)|^2}\Omega(a)(v+\mathcal{Q}_w(E(w,\infty))E(v,\infty)\\&-\frac{(1+\Omega(a)(w+\mathcal{Q}(E(w,\infty)))^2}{|1-i\zeta\partial_\zeta E(w,\infty)|^4}\\&\cdot\Im\left(\overline{(1-i\zeta \partial_\zeta E(w,\infty))}\zeta \partial_\zeta E(v,\infty) \right)+Gv.
\end{aligned}
\end{equation*}
 This resembles the Fréchet derivative of the operator for the infinite depth case. See the discussion following in \cite[ Theorem 3]{paper2}. In fact, we have that $D_w\mathcal{G}[\Im (w(a)),0,a,0]$ coincides with $D_w\mathcal{F}[w(a),0,a]$ written in terms of the imaginary part of $w(a)$ restricted to the unit circle. Here $\mathcal{F}$ denotes the operator defined in \cite[Section 4]{paper2}. This implies that $D_w\mathcal{G}[\Im(w(a)),0,a,0]$ is an isomorphism. By the Implicit Function Theorem, \cite[Theorem 3]{paper2} can be extended in the following way.
\begin{theorem}
\label{imp}
For each $a_0\in(0,\tfrac{1}{4})$ there exists $\epsilon>0$ and a continuous operator
\begin{equation*}
W:(-\epsilon,\epsilon)\times(a_0-\epsilon,a_0+\epsilon)\times (-\epsilon,\epsilon)\rightarrow U
\end{equation*}
 such that $W(0,a,0)=\Im (w(a))$ and 
\begin{equation*}
\mathcal{G}(W(G,a,l);G,a,l)=0.
\end{equation*}
 Moreover, there exists a $\delta>0$ such that for all $(G,a,l)\in (-\epsilon,\epsilon)\times(a_0-\epsilon,a_0+\epsilon)\times (-\epsilon,\epsilon)$, the following statements are equivalent.
\renewcommand{\theenumi}{(\roman{enumi})}
\begin{enumerate}
    \item $\mathcal{G}(w;G,a,l)=0$ and $||w-\Im (w(a))||_X<\delta$;
    \item $w=W(G,a,l)$.
\end{enumerate}
\end{theorem}
By the same argument as in \cite{paper2}, Theorem \ref{imp} implies the following two results.
\begin{theorem}[Overhanging waves]
For all $a\in(a_{\text{crit}},a_{\text{max}})$, there exists an $\epsilon>0$, such that for all $G\in(-\epsilon,\epsilon)$ and all $H$ sufficiently large, there exists a solution of (\ref{conditions}) and (\ref{kinematic at the finite bottom}) with $\Omega=\Omega(a)$ and $B=B(a)$ whose profile is overhanging.
\end{theorem}
\begin{theorem}[Touching waves]
There exists an $\epsilon>0$, such that for all $G\in(-\epsilon,\epsilon)$ and all $H$ sufficiently large, there exists a solution of (\ref{conditions}) and (\ref{kinematic at the finite bottom}) whose profile intersects itself tangentially, enclosing a small bubble of air.
\end{theorem}
We note that these results were mentioned in Section 5 of \cite{paper2}.
\section{Continuous curve of solutions}\label{section: cont curve}
\par This section aims to demonstrate the existence of a continuous curve of gravity waves connecting laminar flow and touching waves. We start by stating our results in the first subsection. Subsequently, in the second subsection, we establish a uniform version of Theorem \ref{imp}.
\subsection{Statement of results}
\par We will prove a uniform version of Theorem \ref{imp} in the second subsection, but before that, let us state the theorem.
\begin{theorem}[Main Theorem]\label{impuni finite}
\par For all $\gamma\in(0,\tfrac{1}{8})$ there exists $\epsilon>0$ and a continuous operator
\begin{equation*}
W:S\rightarrow U,
\end{equation*}
 where
\begin{equation*}
S:=\bigg\{(G,a,l)\in(-\epsilon,\epsilon)\times(-\epsilon,\tfrac{1}{4}-\gamma]\times(-\epsilon,\epsilon): G\geq \Tilde{G}(a,l)\bigg\}.
\end{equation*}
 Here, $\Tilde{G}$ is a continuous surjective map between a two-dimensional and a one-dimensional neighbourhood of the origin satisfying $\Tilde{G}(0,0)=0$. The operator satisfies the following three conditions.
\begin{enumerate}
    \item $W(0,a,l)=\Im(w(a))$;
    \item $W(\Tilde{G}(a,l),a,l)$ is a constant function of $\zeta$;
    \item $\mathcal{G}(W(G,a,l);G,a,l)=0$.
\end{enumerate}

 Moreover, there exists a $\delta>0$ such that for all $(G,a,l)\in S$, the following statements are equivalent:
\renewcommand{\theenumi}{(\roman{enumi})}
\begin{enumerate}
    \item $\mathcal{G}(w;G,a,l)=0$ and $||w-\Im(w(a))||_X<\delta$;
    \item $w=W(G,a,l)$.
\end{enumerate}
\end{theorem}
\par The solutions from Theorem \ref{impuni finite} give rise to solutions of \eqref{conditions pre} and \eqref{bottom conditions pre}, allowing us to state the following results. 
\begin{theorem}\label{continuous curve}
There exists an $\epsilon>0$ such that for all $G\in(-\epsilon,\epsilon)$ fixed, there exists a continuous curve of solutions of (\ref{conditions}) and (\ref{kinematic at the infinity bottom}) between a laminar flow and a touching wave.
\end{theorem}
\begin{corollary}[Breaking Waves]\label{vertical wave}
For every sufficiently small $\epsilon>0$  and $G\in(-\epsilon,\epsilon)$ fixed, there exists a solution of (\ref{conditions}) and (\ref{kinematic at the infinity bottom}) whose profile is vertical at a point but never overhanging.
\end{corollary}
\begin{theorem}\label{continuous curve finite}
There exists an $\epsilon>0$ such that for all $G\in(-\epsilon,\epsilon)$ fixed and all $H$ sufficiently large, there exists a continuous curve of solutions of (\ref{conditions}) and (\ref{kinematic at the finite bottom}) between a laminar flow and a touching wave.
\end{theorem}
\begin{corollary}[Finite depth breaking waves]\label{vertical wave finite}
For every sufficiently small $\epsilon>0$, $G\in(-\epsilon,\epsilon)$ and all $H$ sufficiently large, there exists a solution of (\ref{conditions}) and (\ref{kinematic at the finite bottom}) whose profile is vertical at a point but never overhanging.
\end{corollary}
\subsection{Proof of Main Theorem}
We start by using compactness to establish a version of Theorem \ref{impuni finite} for compact subsets of $(0,\tfrac{1}{4})$. Then, we compute constant solutions of the operator problem. Subsequently, we find a curve of parameters such that the Fréchet Derivative of the operator at every point along the curve has the same kernel as at the point $(0,0,0)$. Next, we use a local bifurcation result to construct solutions near the point $(0,0,0)$. Finally, we connect this curve of solutions with the solutions obtained using compactness.

\begin{lemma}\label{lemmaimp finite}
\par For all $\lambda\in(0,\tfrac{1}{8})$ there exist $\epsilon>0$ and a real-analytic operator
\begin{equation*}
W:(-\epsilon,\epsilon)\times[\lambda,\tfrac{1}{4}-\lambda]\times(-\epsilon,\epsilon)\rightarrow U
\end{equation*}
 such that $W(0,a,0)=\Im (w(a))$ and 
\begin{equation*}
\mathcal{G}(W(G,a,l);G,a,l)=0.
\end{equation*}
 Moreover, there exists a $\delta>0$ such that for all $(G,a,l)\in (-\epsilon,\epsilon)\times[\lambda,\tfrac{1}{4}-\lambda]\times(-\epsilon,\epsilon)$ the following statements are equivalent.
\renewcommand{\theenumi}{(\roman{enumi})}
\begin{enumerate}
    \item $\mathcal{G}(w;G,a,l)=0$ and $||w-\Im (w(a))||_X<\delta$;
    \item $w=W(G,a,l)$.
\end{enumerate}
\end{lemma}
\begin{proof}
\par We start by applying Theorem \ref{imp} to each $a\in[\lambda, \tfrac{1}{4}-\lambda]$, yielding an $\epsilon_a$ and a $\delta_a$ satisfying the conditions in the theorem. Afterwards, we cover $[\lambda,\tfrac{1}{4}-\lambda]$ with $\bigcup_{a\in[\lambda,\tfrac{1}{4}-\lambda]}\bigg(a-\epsilon_a, a+\epsilon_a\bigg)$. By the Heine-Borel Theorem, we can take a finite subcover, meaning
\begin{equation*}
\big[\lambda, \tfrac{1}{4}-\lambda\big]\subset\bigcup_{k=1}^n\bigg(a_k-\epsilon_{a_k}, a_k+\epsilon_{a_k}\bigg).
\end{equation*}
Afterwards, we define $\delta:=\min_{k=1,..,n}\delta_{a_k}$. Let
\begin{equation}\label{operator wk}
W^k:(-\epsilon_{a_k},\epsilon_{a_k})\times (a_k-\epsilon_{a_k},a_k+\epsilon_{a_k})\times(-\epsilon_{a_k},\epsilon_{a_k})\rightarrow U
\end{equation}
be the operator obtained from applying Theorem \ref{imp} to $a_k$. This operator satisfies
\begin{equation*}
||W^k(G,a,l)-\Im(w(a))||_X<\delta_{a_k},\qquad\text{for }|a-a_k|<\epsilon_{a_k},\quad|G|<\epsilon_{a_k}\text{  and  } |l|<\epsilon_{a_k}.
\end{equation*}
 Since $W^k$ is continuous, we can choose $\Tilde{\epsilon}_{a_k}$ such that
\begin{equation*}
    ||W^k(G,a,l)-\Im(w(a))||_X<\delta,\qquad\text{for }|a-a_k|<\epsilon_{a_k},\quad|G|<\Tilde{\epsilon}_{a_k}\text{  and  } |l|<\Tilde{\epsilon}_{a_k}.
\end{equation*}
Then, we define
\begin{equation*}
\epsilon:=\min_{k=1,..,n}\Tilde{\epsilon}_{a_k}.
\end{equation*}
Next, we define the operator $W:(-\epsilon,\epsilon)\times[\lambda,\tfrac{1}{4}-\lambda]\times (-\epsilon,\epsilon)\rightarrow U$ by
\begin{equation}\label{wdefinition}
W(G,a,l):= W^k(G,a,l) \qquad \text{if } |a-a_k|<\epsilon_{a_k}, k=1,...,n.
\end{equation}
 Let $a$ satisfy $|a-a_j|<\epsilon_{a_j}$ and $|a-a_k|<\epsilon_{a_k}$ for some $j$ and $k$, $G$ satisfy $|G|<\epsilon$ and $l$ satisfy $|l|<\epsilon$. Then $(G,a,l)$ is in the domain of both $W^j$ and $W^k$. Since
\begin{equation*}
||W^i(G,a,l)-\Im(w(a))||_X<\delta\qquad\text{for }i=j,k,
\end{equation*}
it follows from the uniqueness of Theorem \ref{imp} that $W^j(G,a,l)=W^k(G,a,l)$, and therefore $W$ is well defined.
\end{proof}
\par We turn our attention now to the bifurcation point, starting by computing the partial Fréchet derivative at $(w(0),0,0,0)$,
\begin{equation}\label{frech at zero}
D_w\mathcal{G}[\Im(w(0));0,0,0]v=\Im (E(v,\infty)-\zeta E(v,\infty)_\zeta)=-\Im\Big(\sum_{n=0}^\infty (n-1)v_n\zeta^n\Big),
\end{equation}
 which is not invertible. Note that the symmetry assumption imposed on $X$ implies the coefficients $v_n$ are purely imaginary. In fact, it follows from (\ref{frech at zero}) that the kernel and the complement of the range of $D_w\mathcal{G}[w(0);0,0,0]$ are both one-dimensional. The kernel is spanned by $i\zeta$, and the complement of its range is spanned by $\Re(\zeta)$.

\par The next step is to look for constant solutions of $\mathcal{G}(w;G,a,l)=0$, which correspond to laminar flow solutions of (\ref{conditions}).
\begin{lemma}\label{lemmaconst}
There exists a real-analytic operator $d:\mathbb{R}\times(-\infty,1/3)\rightarrow \{iR:R\in\mathbb{R}\}$ such that
\begin{equation*}
\mathcal{G}(d(G,a);G,a,l)=0\qquad\qquad\text{and}\qquad\qquad d(0,0)=0
\end{equation*}
 for $(G,a)\in\mathbb{R}\times(-\infty,1/3)$ and $l\in\mathbb{R}$. Moreover, there exist open neighbourhoods $V\subset \mathbb{R}^2$, $T\subset\mathbb{R}$ containing zero and a smooth map $\Tilde{G}:V\rightarrow T$, such that
\begin{equation*}
\ker D_x\mathcal{G}\Big[d(\Tilde{G}(a,l),a); \Tilde{G}(a,l), a,l\Big]=\langle i\zeta\rangle_X
\end{equation*}
 for $(a,l)\in V$, $\Tilde{G}(0,0)=0$ and for every fixed $l$ the map $\Tilde{G}_l(a):= \Tilde{G}(a,l)$ is invertible.
\end{lemma}
\begin{proof}
\par Let $c$ be a fixed constant. A computation reveals
\begin{equation*}
\mathcal{G}(ic;G,a,l)=\frac{1}{2}(1+\Omega(a)c)^2+Gc-B(a).
\end{equation*}
 Thus $\mathcal{G}(ic;G,a,l)=0$ if and only if
\begin{equation}\label{d}
c=-\frac{\Omega(a)+G}{\Omega(a)^2}\pm\frac{1}{\Omega(a)^2}\sqrt{G^2+2\Omega(a)G+2 B(a)\Omega(a)^2}.
\end{equation}
 Evaluating this expression at $a=G=0$, we find that when the sign before the square root is positive, $c=0$. We define the operator $d:\mathbb{R}\times(-\infty,1/3)\rightarrow \{iR:R\in\mathbb{R}\}$ by
\begin{equation}\label{d constant d}
d(G,a)=i\bigg(-\frac{\Omega(a)+G}{\Omega(a)^2}+\frac{1}{\Omega(a)^2}\sqrt{G^2+2\Omega(a)G+2 B(a)\Omega(a)^2}\bigg).
\end{equation}
\par We turn our focus to the second part of the lemma and compute the Fréchet derivative of (\ref{more general of operators}) at $(d(G,a); G,a,l)$,
\begin{equation}
\begin{aligned}\label{lin}
D_w\mathcal{G}[d(G,a); G, a,l]v=&\Im\bigg(\Big((1+\Omega(a)d(G,a))\Omega(a)+G\Big)E(v,\frac{1}{l^2})
\\&-(1+\Omega(a)d(G,a))^2\zeta \partial_\zeta E(v,\frac{1}{l^2})\bigg).
\end{aligned}
\end{equation}
 By plugging in $v=i\zeta^k$ in (\ref{lin}) and equating it to zero, we obtain
 \begin{equation*}
 \big((1 + \Omega(a) \Im(d(G,a)))\Omega(a) + G\big)E(v,\frac{1}{l^2})-k(1 + \Omega(a) \Im(d(G,a)))^2 E(v,\frac{1}{l^2})=0,
 \end{equation*}
which implies that $i\zeta$ spans the kernel of $D_w\mathcal{G}[d(G,a); G, a,l]$ if and only if 
\begin{equation}\label{bif}
\big((1 + \Omega(a) \Im(d(G,a)))\Omega(a) + G\big)E(v,\frac{1}{l^2})-(1 + \Omega(a) \Im(d(G,a)))^2 E(v,\frac{1}{l^2})=0.
\end{equation}
 We note that in a neighbourhood of $(0,0,0)$, (\ref{d constant d}) can be described in the following way,
\begin{equation}\label{deltad}
\Im(d(G,a))=4a(1+\mathcal{O}(a+G)).
\end{equation}
 Therefore, we can compute the partial derivative with respect to $G$ of the left-hand side of (\ref{bif}) evaluated at $(G,a,l)=(0,0,0)$ to be equal to $1$. An application of the Implicit Function Theorem yields $\Tilde{G}$ as described in the lemma. Fixing $l$ and applying the Implicit Function Theorem with respect to $a$ will yield the invertibility claim.
\end{proof}

\begin{lemma}\label{bifurcation from all points finite}
\par There exist $\epsilon>0$ and a smooth operator $v:\Tilde{S}\rightarrow U$, where
\begin{equation*}
\Tilde{S}:= \left\{(G,a,l)\in \mathbb{R}\times (-\epsilon,\epsilon)\times  (-\epsilon,\epsilon):\Tilde{G}(a,l)\leq G< \Tilde{G}(a+\epsilon,l+\epsilon)\right\}
\end{equation*}
such that $v(G,a,l)$ satisfies the equation
\begin{equation*}
\mathcal{G}\left(v(G,a,l)+d(G,a);G,a,l\right)=0
\end{equation*}
 for $(G,a,l)\in \Tilde{S}$. Moreover, $v(G,a,l)$ is of the form
\begin{equation}\label{formula for v g a}
v(G,a,l)=C_{G,l}(a-\Tilde{G}_l^{-1}(G))^\frac{1}{2}\Re\zeta+\mathcal{O}((a-\Tilde{G}_l^{-1}(G))+l),
\end{equation}
 where $C_{G,l}$ is a positive constant that depends continuously on $G$ and $l$.
Furthermore, there exists a $\delta>0$ such that $v(G,a,l)+d(G,a)$ is the only nontrivial solution of $\mathcal{G}\big(v,G,a,l\big)=0$, satisfying $||v||_X<\delta$.
\end{lemma}
\begin{proof}
We start by defining the auxiliary operator
$\mathcal{K}:U\times\mathbb{R}\times(-\infty,1/3)\times\mathbb{R}\rightarrow Y$ by
\begin{equation}\label{k}
\mathcal{K}(v;G,a,l)=\mathcal{G}( d(G,a)+v;G,a,l).
\end{equation}
It follows that $\mathcal{K}(0;G,a,l)=0$, for all $(G,a,l)\in\mathbb{R}\times(-\infty,\tfrac{1}{3})\times\mathbb{R}$. Moreover, 
\begin{equation*}
\ker D_v\mathcal{K}\left(0;\Tilde{G}(a_{0},l_0),a_0,l_0)\right)=\langle i\zeta\rangle_X\quad\text{and}\quad \ran D_v\mathcal{K}\left(0;\Tilde{G}(a_{0},l_0),a_0,l_0\right)=Y\ominus \langle \Re(\zeta)\rangle_Y
\end{equation*} for all $(a_0,l_0)\in V$. Next, we compute the mixed derivative
\begin{equation}\label{transv}
D^2_{v a}\mathcal{K}\left(0;\Tilde{G}(a_{0},0),a_0,0\right)(i\zeta)=\Big(-2+\mathcal{O}(a_0)\Big)\Re\zeta.
\end{equation}
 Varying $l$ in a neighbourhood of $0$ will only create smooth dependence on $l$ (see (A.2) and the discussion following it in \cite{conststrauvarv}), the same is true for the other partial derivatives below, so without loss of generality, we consider $l=0$. The expression (\ref{transv}) needs to be nonzero in order to satisfy the transversality condition. To achieve this, we redefine $V$ if necessary. Then, we can apply Corollary \ref{crandy 3} to obtain real analytic operators
\begin{equation}\label{v and a proof step}
\Tilde{v}:(-\Tilde{\epsilon},\Tilde{\epsilon})^3\rightarrow U\qquad\text{and}\qquad \Tilde{a}:(-\Tilde{\epsilon},\Tilde{\epsilon})^3\rightarrow\mathbb{R},
\end{equation}
such that 
\begin{equation*}
\mathcal{K}\Big(\Tilde{v}(s,G,l),G,\Tilde{a}(s,G,l),l\Big)=0, \text{for all }|s|,|G|,|l|<\Tilde{\epsilon},
\end{equation*}
 $\Tilde{v}(s,G,l)=si\zeta+\mathcal{O}(s^2+sG+sl)$. We can then use the bifurcation formulas from Section I.6 of \cite{kielhöfer2011bifurcation} to compute the derivatives of $\Tilde{a}(s,G,l)$ with respect to $s$. We proceed with computing the second-order Fréchet derivative at $(0,\Tilde{G}(a_0,0),a_0,0)$,
\begin{align}\label{double v frechet}
D^2_{vv}\mathcal{K}(0,\Tilde{G}(a_0,0),a_0,0)(u)^2= \Omega&(a_0)\mathcal{Q}_{ww}(d(\Tilde{G}(a_0,0),a))E(u,0)^2\nonumber\\  -&4\Omega(a_0)\Im(E(u,0))\Im(\zeta E(u,0)_\zeta)+4\Big(\Im(\zeta E(u,0)_\zeta)\Big)^2\\ -&\Im(i|E(u,0)_\zeta|^2),\nonumber
\end{align}
 where we use $\mathcal{Q}_{w}(d(G,a))(u)=0$. Evaluating (\ref{double v frechet}) at $u=i\zeta$ reveals
\begin{equation}\label{second fréchet v}
D^2_{vv}\mathcal{K}(0,\Tilde{G}(a_0),a_0,0)(i\zeta)^2=\Omega(a_0)-4\Omega(a_0)\Re(\zeta)^2+4\Re(\zeta)^2-1,
\end{equation}
 where we used $\mathcal{Q}_{w}(d(G,a))(v)=0$ and $\mathcal{Q}_{ww}(d(G,a))(i\zeta)^2=1$. Note that (\ref{second fréchet v}) has zero projection on $\langle \Re(\zeta)\rangle_Y$ and therefore, $\tfrac{\partial}{\partial s}\Tilde{a}(0,\Tilde{G}(a_0))=0$. Subsequently, we compute the third-order linearization in $v$ at $(0,\Tilde{G}(a_0),a_0)$,
\begin{align}\label{triple v Frechet}
D^3_{vvv}\mathcal{K}(0,\Tilde{G}(a_0),a_0,0)(u)^3&=3\Omega(a_0)^2(\Im u)\mathcal{Q}_{ww}(u)^2\nonumber\\&-6\Omega(a_0)\mathcal{Q}_{ww}(u)^2\Im(\zeta u_\zeta)-6\Omega(a_0)^2(\Im u)^2\Im(\zeta u_\zeta)\\&+3\Omega(a_0)\Im(u)\bigg(8\Big(\Im(\zeta u_\zeta)\Big)^2-2\Im(i|u_\zeta|^2)\bigg)\\&-24\Big(\Im(\zeta u_\zeta)\Big)^3+12\Im(\zeta u_\zeta)\Im(i|u_\zeta|^2).\nonumber
\end{align}
 Evaluating (\ref{triple v Frechet}) at $u=i\zeta$ reveals
\begin{equation*}
D^3_{vvv}\mathcal{K}(0,\Tilde{G}(a_0),a_0,0)(i\zeta)^3=3\big(\Omega(a_0)-2\big)^2\Re\zeta-6(\Omega(a_0)-2)^2(\Re\zeta)^3,
\end{equation*}
which has nonzero projection on $\langle \Re(\zeta)\rangle_Y$. It follows that $\tfrac{\partial^2}{\partial s^2}\Tilde{a}(0,\Tilde{G}(a_0),0)=\tfrac{(\Omega(a_0)-2)^2}{2+\mathcal{O}(a_0)}$. We, again, possibly, redefine $V$ to ensure $\tfrac{\partial^2}{\partial s^2}\Tilde{a}(0,\Tilde{G}(a_0),0)>0$, for all $a_0\in V$. Since $\tfrac{\partial^2}{\partial s^2}\Tilde{a}(0,\Tilde{G}(a_0),0)>0$, it follows that we can choose an $\epsilon$ such that the real analytic operator $\Tilde{a}(\cdot,\Tilde{G}(a_0),0)$ is injective on the set $[0,\epsilon)$. We define $v:(-\epsilon,\epsilon)\times[a_0,a_0+\epsilon)\times (-\epsilon,\epsilon)\rightarrow U$ by
\begin{equation*}
v(G,a,l):= \Tilde{v}\left(\Tilde{a}^{-1}(s,G,l)\right).
\end{equation*}
 The formula in (\ref{formula for v g a}) follows directly from the definition above.
\end{proof}

\begin{proof}[Proof of Theorem \ref{impuni finite}]
\par Let $\gamma>0$ be given. Let $v,\epsilon_1,\delta_1$ and $\Tilde{S}$ be as described in Lemma \ref{bifurcation from all points finite}. Next, we note that for small and positive $a$, (\ref{sol of F 0 grav}) can be expressed as
\begin{equation}
w(a)=-\frac{4i\sqrt{a}\zeta}{1+\sqrt{a}\zeta}=-4i\sqrt{a}\zeta+\mathcal{O}(a).
\end{equation}
 Therefore, we can choose a positive and small enough $\lambda$ such that $||w(\lambda)||_X<\tfrac{\delta_1}{2}$, $||d(0,\lambda)||_X<\tfrac{\delta_1}{2}$, $\lambda<\epsilon_1$ and $\lambda<\gamma$. We then apply Lemma \ref{lemmaimp finite} for such a $\lambda$, which yields $\epsilon_2,\delta_2$ and $W^{\lambda}$ as described in the lemma. We then define
\begin{equation*}
\delta:=\min(\delta_1,\delta_2).
\end{equation*}
We choose $\Tilde{\epsilon}_1$ and $\Tilde{\epsilon}_2$ to be small enough such that
\begin{equation*}
||v(G,a,l)-d(G,a)||_X<\delta \qquad\text{ for }|G|,|l|<\Tilde{\epsilon}_1
\end{equation*}
 and
\begin{equation*}
||W^\lambda(G,a)-w(a)||_X<\delta \qquad\text{ for }|G|,|l|<\Tilde{\epsilon}_2.
\end{equation*}
 We then define
\begin{equation*}
\epsilon:= \min(\Tilde{\epsilon}_1,\Tilde{\epsilon}_2)
\end{equation*}
and
\begin{equation*}
S:=\bigg\{(G,a)\in(-\epsilon,\epsilon)\times(-\infty,\tfrac{1}{\sqrt{2}}-\gamma]:a\geq \Tilde{G}^{-1}(G)\bigg\}.
\end{equation*}
Let $W:S\rightarrow U$ be defined by
\begin{equation}\label{final w}
W(G,a):= \begin{cases}
  W^\lambda(G,a)  & \text{if } a>\lambda\\
  v(G,a)+d(G,a) & \text{if }(G,a)\in\Tilde{S}\cap S.
\end{cases}
\end{equation}
Next, we show (\ref{final w}) is well-defined. Let $(G,a)\in \Tilde{S}$ and $a>\lambda$, then 
\begin{equation*}
\mathcal{K}\Big(W_\lambda(G,a)-d(G,a),G,a\Big)=0    
\end{equation*}
 and
\begin{equation*}
||W_\lambda(G,a)-d(G,a)||_X<\frac{\delta(\mu)}{2}+\frac{\delta(\mu)}{2}\leq \delta(\mu).
\end{equation*}
 By Lemma \ref{bifurcation from all points finite}, the solutions coincide, meaning (\ref{final w}) is well-defined.
\end{proof}
\section{Critical layers}\label{section: crit layers}
Critical layers are defined as the set of points where the horizontal velocity,
\begin{equation}\label{psi_y def crit layer}
\psi_y=\frac{\psi_\alpha y_\alpha +\psi_\beta x_\alpha}{|z_\alpha|^2}    
\end{equation}
vanishes. Our next result characterizes the local behaviour of a critical layer when it intersects at a point on the surface where the tangent is vertical but which is not a stagnation point. Notably, the results hold for a general water wave and not just the ones we constructed.
\begin{proposition}\label{critical layer prop}
Let $z$ denote the conformal map to the half plane of a solution to (\ref{conditions}) with $\tau=0$ and assume the same conditions as in Proposition \ref{analytic proposition}. Additionally, assume that $\alpha_{\text{crit}}$ is a point in the $(\alpha,\beta)$ plane at which the wave profile is vertical.
\renewcommand{\theenumi}{(\roman{enumi})}
\begin{enumerate}
\item If the function $x(\alpha+i0)$ does not have a local extremum at $\alpha_{\text{crit}}$ and $G>0$, then no critical layers touch the free surface at $\alpha_{\text{crit}}$. However, their analytic extension from outside the fluid domain do touch $\alpha_{\text{crit}}$;
\item If the function $x(\alpha+i0)$ does not have a local extremum at $\alpha_{\text{crit}}$ and $G<0$, there exist critical layers from the water touching the point $\alpha_{\text{crit}}$. Moreover, their analytic extension from outside the fluid domain do not touch $\alpha_{\text{crit}}$;
\item If the function $x(\alpha+i0)$ has a local extremum at $\alpha_{\text{crit}}$ and $G\neq 0$, there exist critical layers from the water touching the point $\alpha_{\text{crit}}$. Moreover, their analytic extension from outside the fluid domain also touch $\alpha_{\text{crit}}$.
 \end{enumerate}
\end{proposition}
\begin{proof}
Without loss of generality, assume the wave is travelling in the positive $x$ direction. Let $\alpha_{\text{crit}}$ denote the point at which the wave is vertical. It follows that
\begin{equation*}
x_\alpha(\alpha_{\text{crit}},0)=0
\end{equation*}
and the smallest positive integer $k$, such that $\frac{\partial^k}{\partial\alpha^k}x(\alpha_{\text{crit}},0)$ does not vanish, must be even if the wave is a breaking wave and odd if the wave overhangs at the point. Moreover
\begin{equation}\label{derivative cases}
\frac{\partial^k}{\partial\alpha^k}x(\alpha_{\text{crit}},0)
>0
\end{equation}
 if $k$ is odd. It follows from (\ref{psi_y def crit layer}) that 
\begin{equation}\label{smth}
\psi_y=0\iff \psi_\alpha y_\alpha+\psi_\beta x_\alpha=0.
\end{equation}
Let $U$ be the open set from Proposition \ref{analytic proposition} where $\psi$ and $z$ have a real-analytic extension and define $F:U\rightarrow\mathbb{R}$ by
\begin{equation}\label{fff}
F(\alpha,\beta)=\psi_\alpha(\alpha,\beta) y_\alpha(\alpha,\beta)+\psi_\beta (\alpha,\beta)x_\alpha(\alpha,\beta).
\end{equation}
(\ref{smth}) implies $(\alpha,\beta)$ belongs to the critical layer if and only if $F(\alpha,\beta)=0$. Since $F(\alpha_\text{crit},0)=0$, we want to use Taylor expansions of (\ref{fff}) near the point $(\alpha_\text{crit},0)$ to figure out the local behaviour of the critical layers.

 Notice that $\psi$ is constant on $\{\beta=0\}$, and therefore
\begin{equation}\label{no more alphas}
\frac{\partial^j}{\partial\alpha^j}\psi(\alpha_\text{crit},0)=0,\text{ for all }j\in \mathbb{Z}^+.
\end{equation}
 Next, we rewrite (\ref{dynamic boundary condition in conditions}) as
\begin{equation}\label{berncrit}
\psi_\beta=\pm \sqrt{2(x_\alpha^2+y_\alpha^2)(B-G y)}\qquad\text{on }\{\beta=0\}
\end{equation}
 for later use. Given all the $\alpha$ derivatives of $\psi$ vanish at the critical point, exploring alternative branches of (\ref{berncrit}) results in the partial derivatives of $F$ at the critical point differing only in sign. Consequently, we can consider the positive branch of equation (\ref{berncrit}) without any loss of generality.

\par We proceed to compute the partial derivatives of $F$ at the point $(\alpha_\text{crit},0)$. Starting by
\begin{equation*}
\begin{split}
\frac{\partial^j}{\partial\alpha^j}F(\alpha,\beta)=\sum_{l=0}^j\binom{j}{l}\left(\frac{\partial^l\psi_\alpha}{\partial\alpha^l}\frac{\partial^{j-l}y_\alpha}{\partial\alpha^{j-l}}+\frac{\partial^l\psi_\beta}{\partial\alpha^l}\frac{\partial^{j-l}x_\alpha}{\partial\alpha^{j-l}}\right),\text{ for all }j\in \mathbb{Z}^+.
\end{split}
\end{equation*}
 Using (\ref{no more alphas}) and $\frac{\partial^j}{\partial\alpha^j}x=0$ for $j<k$ to evaluate the expression above at $(\alpha_{\text{crit}},0)$ yields
\begin{equation}\label{alpha does not disappear}
\frac{\partial^j}{\partial\alpha^j}F(\alpha_{\text{crit}},0)=\begin{cases}0,&\text{if }j<k-1, \\
\psi_\beta\frac{\partial^kx}{\partial\alpha^k},&\text{if }j=k-1.
\end{cases}
\end{equation}

 Next, we compute the derivatives with respect to $\beta$
\begin{equation}\label{beta derivatives of F}
\begin{aligned}
\frac{\partial^j}{\partial\beta^j}F(\alpha,\beta)&=\sum_{l=0}^j\binom{j}{l}\left(\frac{\partial^l\psi_\alpha}{\partial\beta^l}\frac{\partial^{j-l}y_\alpha}{\partial\beta^{j-l}}+\frac{\partial^l\psi_\beta}{\partial\beta^l}\frac{\partial^{j-l}x_\alpha}{\partial\beta^{j-l}}\right).\\
\end{aligned}
\end{equation}
 Plugging in $(\alpha_{\text{crit}},0)$ for $j=1$ above yields
\begin{equation}\label{beta only with g}
\begin{aligned}
\frac{\partial}{\partial\beta}F(\alpha_\text{crit},0)&=\psi_{\alpha\beta}y_\alpha-\psi_\beta y_{\alpha\alpha}\\&=\left(\sqrt{2(B-G y)}y_{\alpha\alpha}-\frac{G y_\alpha^2 }{\sqrt{2(B-G y)}}\right)y_\alpha-\sqrt{2(B-G y)}y_\alpha y_{\alpha\alpha}\\&=-\frac{G y_\alpha^3 }{\sqrt{2(B-G y)}}.
\end{aligned}
\end{equation}
 We utilized equation (\ref{berncrit}) above to calculate the derivatives of $\psi_\beta$. Additionally, we assume $\sqrt{y_\alpha^2}=y_\alpha$ because symmetry implies the existence of a critical point where $y_\alpha$ is positive and another where it is negative. The Taylor expansion of $F$ at the point $(\alpha_\text{crit},0)$ is then
\begin{equation}\label{main taylor}
\begin{aligned}
&\left(-\frac{G y_\alpha^3 }{\sqrt{2(B-G y)}}\right)\Bigg|_{\alpha=\alpha_{\text{crit}}}\beta+\left(\sqrt{2(B-G y)} y_\alpha\frac{\partial^kx}{\partial \alpha^k}\right)\Bigg|_{\alpha=\alpha_{\text{crit}}}(\alpha-\alpha_{\text{crit}})^{k-1}\\=& \mathcal{O}\Big(|\alpha-\alpha_{\text{crit}}|^k+|\beta|^2+|(\alpha-\alpha_{\text{crit}})\beta|\Big).
\end{aligned}
\end{equation}
 When $k$ is odd, (\ref{main taylor}) has the form
\begin{equation}\label{breakingtaylorexpansion}
-C_1 G \beta+C_2(\alpha-\alpha_{\text{crit}})^{k-1}=\mathcal{O}\Big(|\alpha-\alpha_{\text{crit}}|^k+|\beta|^2+|(\alpha-\alpha_{\text{crit}})\beta|\Big),
\end{equation}
 where $C_1$ and $C_2$ denote positive constants that depend on $G$. It follows from applying the Implicit Function Theorem to (\ref{breakingtaylorexpansion}) that solutions of $F(\alpha,\beta)=0$ in a neighbourhood of $(\alpha_\text{crit},0)$ have the form 
\begin{equation}
\beta=\frac{C_2}{C_1 G} (\alpha-\alpha_\text{crit})^{k-1} + \mathcal{O}( |\alpha-\alpha_\text{crit}|^{k}).
\end{equation}
 Since the exponent $k-1$ is even, statements $(i)$ and $(ii)$ follow.
 When $k$ is even, (\ref{main taylor}) has the form
\begin{equation}\label{overhangingtaylorexpansion}
C_3 G \beta+C_4(\alpha-\alpha_{\text{crit}})^{k-1}=\mathcal{O}\Big(|\alpha-\alpha_{\text{crit}}|^k+|\beta|^2+|(\alpha-\alpha_{\text{crit}})\beta|\Big),
\end{equation}
 where $C_3$ and $C_4$ denote constants that depend on $G$. It follows from applying the Implicit Function Theorem to (\ref{overhangingtaylorexpansion}) that solutions of $F(\alpha,\beta)=0$ in a neighbourhood of $(\alpha_\text{crit},0)$ have the form 
\begin{equation}
\beta=\frac{C_4}{C_3 G} (\alpha-\alpha_\text{crit})^{k-1} + \mathcal{O}( |\alpha-\alpha_\text{crit}|^{k}).
\end{equation} Since the exponent $k-1$ is odd, statement $(iii)$ follows.
\end{proof}

\section{Local Bifurcation Theory}
In this appendix, we provide the local bifurcation theory used throughout the paper. The proposition below is a modified version of the Crandall--Rabinowitz Theorem \cite[Theorem 8.3.1]{10.2307/j.ctt1dwstqb} with two additional parameters. The proof is very similar to that in the reference.
\begin{proposition}\label{crandy two}
Let $X$ and $Y$ be Banach spaces, $\mu_0,\lambda_0,\nu_0\in \mathbb{R}$ and $F:X\times \mathbb{R}\times\mathbb{R}\times\mathbb{R}\rightarrow Y$ an analytic operator satisfying the following three properties.
\renewcommand{\theenumi}{(\roman{enumi})}
\begin{enumerate}
    \item $F(0,\mu,\lambda,\nu)=0$, for all $\mu,\lambda,\nu\in \mathbb{R}$;
    \item $D_x F[0, \mu_0,\lambda_0,\nu_0]$ is a Fredholm operator of index zero with a one-dimensional kernel spanned by $x_0\in X$;
    \item $D_{x \mu}F[0,\mu_0,\lambda_0,\nu_0](x_0)\not \in \ran(D_x F[0, \mu_0,\lambda_0,\nu_0])$ (transversality condition).
\end{enumerate}
\par Then there exist $\epsilon>0$ and a curve $(\chi(s,\lambda,\nu),\mu(s,\lambda,\nu),\lambda,\nu)$, where \begin{align*}
\chi(s,\lambda,\nu)&=s x_0+\mathcal{O}\left(s^2+s(\lambda-\lambda_0)+s(\nu-\nu_0)\right)\\\mu(s,\lambda,\nu)&=\mu_0 +\mathcal{O}\left(s+(\lambda-\lambda_0)+(\nu-\nu_0)\right),\end{align*}such that $F(\chi(s,\lambda,\nu);\mu(s,\lambda,\nu),\lambda,\nu)=0$ for $s,(\lambda-\lambda_0),(\nu-\nu_0)\in (-\epsilon,\epsilon)$. In addition, there exists a neighbourhood of $(0,\mu_0,\lambda_0,\nu_0)$ where these are the only non-trivial solutions.
\end{proposition}

\begin{corollary}\label{crandy 3}
Assuming the same conditions as stated in Proposition \ref{crandy two}, and further, that there exists a continuous curve of parameters $S:=\{(\mu(t),\lambda(t),\nu(t)):t\in(-\delta,\delta)\}$, satisfying the following four properties.
\renewcommand{\theenumi}{(\roman{enumi})}
\begin{enumerate}
    \item $(\mu(0),\lambda(0),\nu(0))=(\mu_0,\lambda_0,\nu_0)$;
    \item For all $t\in(-\delta,\delta)$, $D_x F[0, \mu(t),\lambda(t),\nu(t)]$ is a Fredholm operator of index zero with one-dimensional kernel spanned by $x_0\in X$;
    \item For all $t\in(-\delta,\delta)$, $D_{x \mu}F[0,\mu(t),\lambda(t),\nu(t)](x_0)\not \in \ran(D_x F[0, \mu_0,\lambda_0,\nu_0])$ (transversality condition)
    \item The mapping $t\mapsto D_{x \mu}F[0,\mu(t),\lambda(t),\nu(t)](x_0)$ is continuous.
\end{enumerate}
\par Then, there exists  $\epsilon>0$ and a curve $(\chi(s,\lambda,\nu),\mu(s,\lambda,\nu),\lambda,\nu)$, where \begin{align*}
\chi(s,\lambda,\nu)&=s x_0+\mathcal{O}\left(s^2+s(\lambda-\lambda_0)+s(\nu-\nu_0)\right),\\\mu(s,\lambda,\nu)&=\mu_0 +\mathcal{O}\left(s+(\lambda-\lambda_0)+(\nu-\nu_0)\right)\end{align*}such that $F(\chi(s,\lambda,\nu);\mu(s,\lambda,\nu),\lambda,\nu)=0$ for $s,(\lambda-\lambda_0),(\nu-\nu_0)\in (-\epsilon,\epsilon)$. 

Moreover, for every fixed $(\lambda(t),\nu(t))\in (\lambda_0-\epsilon,\lambda_0+\epsilon)\times (\nu_0-\epsilon,\nu_0+\epsilon)$, the curve \begin{equation*}
(\chi(s,\lambda(t),\nu(t)),\mu(s,\lambda(t),\nu(t))\end{equation*} is the same as the curve obtained by applying the Crandall--Rabinowitz Theorem at the point \begin{equation*}
(0,\mu(t),\lambda(t),\nu(t)).\end{equation*} In addition, there exists a neighbourhood of $(0,\mu_0,\lambda_0,\nu_0)$ where these are the only non-trivial solutions.
\end{corollary}
\par The statement above follows from applying the Proposition \ref{crandy two} at every point along the parameter curve while keeping $\lambda$ fixed, and then shrinking the neighbourhoods so that the uniqueness from Proposition \ref{crandy two} implies that the new solutions are the same as those obtained from Proposition \ref{crandy two}. 
\bibliographystyle{plain}
\bibliography{Gravity_water_waves_over_constant_vorticity_flows}
\end{document}